\documentclass[final, 11pt]{article}
\usepackage[ruled, linesnumbered]{algorithm2e}
\usepackage{mathtools}
\usepackage{mdwlist}
\usepackage{graphicx}
\usepackage{amsthm}
\usepackage{amsfonts}
\usepackage{amsmath}
\usepackage{amssymb}
\usepackage{fancyhdr}
\usepackage{titlesec}
\usepackage{indentfirst}
\usepackage{booktabs}
\usepackage{verbatim}
\usepackage{color}
\usepackage{latexsym}
\usepackage{amscd}
\usepackage{esvect}
\usepackage{graphicx}
\usepackage{upgreek}
\usepackage{subcaption}   
\usepackage{caption}
\usepackage[page,header]{appendix}
\usepackage{titletoc}
\usepackage{mathrsfs}
\usepackage[numbers]{natbib}
\usepackage{float}
\usepackage{lipsum}
\usepackage[export]{adjustbox} 
\numberwithin{equation}{section}
\newtheorem{theorem}{Theorem}[section]

\newtheorem{proposition}[theorem]{Proposition}

\newtheorem{remark}[theorem]{Remark}

\allowdisplaybreaks



\def\e{\varepsilon}

\SetKwInOut{Offline}{Offline}
\SetKwInOut{Online}{Online}
\newtheorem{assumption}{Assumption}

\begin{document}

\title{Error estimate of a bi-fidelity method for kinetic equations with random parameters and multiple scales 
\thanks{The first and the third author is supported by the DOE funding--Simulation Center for Runaway Electron Avoidance and Mitigation. The second author was supported by NSFC grants No. 11871297 and No. 31571071.}
}
\author{Irene M. Gamba\footnote{Oden Institute for Computational Engineering and Sciences and Department of Mathematics, University of Texas at Austin, Austin, TX, USA (gamba@math.utexas.edu)}, Shi Jin\footnote{School of Mathematical Sciences, Institute of Natural Sciences, MOE-LSEC and SHL-MAC, Shanghai Jiao Tong University, Shanghai, China (sjin@wisc.edu)}, Liu Liu\footnote{Oden Institute for Computational Engineering and Sciences and Department of Mathematics, University of Texas at Austin, Austin, TX, USA (lliu@ices.utexas.edu)}}
\date{}
\maketitle

\abstract{In this paper, we conduct uniform error estimates of the bi-fidelity method for multi-scale kinetic equations. We take the Boltzmann and the linear transport equations as important examples, {\color{blue}and discuss various choices of low-fidelity models for more general kinetic equations.} The main analytic tool is the hypocoercivity analysis for kinetic equations, considering solutions in a perturbative setting close to the global equilibrium. This allows us to obtain the error estimates in both kinetic and hydrodynamic regimes. 
}

\section{Introduction}

It is known that solving deterministic kinetic equations are time expensive due to its high-dimensional nature in the phase space. 
When one considers random uncertainties, which are usually described by kinetic equations with random parameters in collision kernels, scattering coefficients, initial or boundary data, forcing or source terms \cite{dimarco2019multi, DPZ-Review, HJS-Boltz, LL-LB, LL-BP, LJ-UQ, SHJ17}, the computational dimension becomes much higher thus it becomes computationally daunting.
Take the Boltzmann equation with multi-dimensional uncertain variables as an example, although one can use the standard stochastic collocation method, often used to deal with high-dimensional parametric PDEs \cite{multilevel,GWZ14,Webster, XH05}, it is not desirable to numerically solve the Boltzmann equation repeatedly, especially given the fact that the nonlinear, non-local collision operator in an integral form is so complicated. 

Often there exist some approximated or less complex low-fidelity models, which 
contain simplified physics or are calculated on a coarser mesh in the physical space thus have cheaper computational costs.   
Although their accuracy may not be high, the low-fidelity models are designed in such a way that they can capture certain important features of the underlying problem and produce reliable predictions. 
In the fluid regime, we are naturally motivated to choose the corresponding hydrodynamic limit equations of the considered kinetic equation as its low-fidelity 
model. 

In \cite{LZ19}, the authors adapt the bi-fidelity method developed in \cite{NGX14, ZNX14}
to efficiently compute high-fidelity solutions of the Boltzmann equation with multi-dimensional random parameters and multiple scales. 
They take advantage of the multiscale nature of kinetic problem and choose the low-fidelity model as the compressible Euler system and show that by using only $O(1)$ runs of the high-fidelity asymptotic-preserving solver for the Boltzmann equation, 
the bi-fidelity approximation can capture well the macroscopic quantities of the solution to the Boltzmann equation in the multi-dimensional random space. 
{\color{blue} Numerous numerical experiments were shown in \cite{LZ19} to prove the accuracy and efficiency of the proposed bi-fidelity method. }
  
The authors of \cite{LZ19} conduct the accuracy and convergence analysis by splitting the error between the high- and bi-fidelity solutions into two parts, namely
the projection and the remainder. 
By incorporating one's knowledge of regularity of the high-fidelity solution in the random space \cite{DJL,LJ-UQ}, based on the best-$N$ approximation theory and \cite{Cohen15} which provides an upper bound for the Kolmogorov $N$-width, \cite{LZ19} shows that the numerical error between the high- and bi-fidelity solutions decays algebraically with respect to the number of high-fidelity runs $N$. The convergence rate is 
 independent of the dimension of the random space and regularity of the initial data. However, it is known that the multi-fidelity estimate error bounds are rarely sharp. Indeed, as the authors acknowledge in \cite[Remark 4.3]{LZ19}, their error bound derived may not be sharp,
since they adopt the theoretical framework developed in \cite{NGX14} which treats general high- and low-fidelity solutions as vector data instead of including
information on specific models or PDEs. 

{\color{blue}Compared to the analysis result in \cite{LZ19},} we provide a more refined, {\it uniform-in-$\e$} ($\e$ being the Knudsen number) error estimate of the bi-fidelity method for solving a more general class of kinetic equations with high-dimensional uncertainties. 
Our analysis is based on the hypocoercivity analysis for the Boltzmann equation \cite{MB15,Mouhot-Neumann}, later extended to the multiscale Boltzmann equation with uncertainties \cite{DJL2,DJL,LJ-UQ} which provides regularity and long-time behavior of the solution in the random space in the perturbative setting, where the solution is close to the global Maxwellian. 
In the fluid regime, one can easily use the error between the (linearized) Boltzmann equation (the high-fidelity model) and the linearized Euler system, or rather the acoustic equations (the low-fidelity model). 

The significant and most {\it novel} part of our analysis is in the kinetic regime in which the fluid limit is not valid, 
yet, in the perturbative setting, one can still establish an error estimate between the Boltzmann equation and its Euler limit, since {\it the moments of the Boltzmann equation are close to the macroscopic quantities of the Euler equations}. In this regime, one cannot directly adopt the scaling used in \cite{LJ-UQ,Mouhot-Neumann} where the scaling coefficient of the perturbative part and the Knudsen number are the same. We need to use and keep track of a different perturbative parameter ($\delta)$ in our analysis in order to obtain the error in the kinetic regime. This allows us to obtain an error analysis even in the kinetic regime where the fluid limit (the low-fidelity) is in general not a good approximation of the Boltzmann equation. As a consequence we obtain {\it uniform} error estimates, from the kinetic and fluid regimes in the perturbative setting, for the bi-fidelity method. 

Since our analysis only uses the property that the moments of the Boltzmann and its Euler limit are close in the perturbative setting, it also suggests that the bi-fidelity model does not necessarily need to be the fluid limit equation. Any moment models in velocity, for example moment closure models, that possess approximately accurate moments as the Boltzmann in the perturbative setting could be used as the low-fidelity model. This will be studied in \ref{sec:3} where we use the linearized Boltzmann and linear transport equations as examples, and also treated the low-fidelity model simply to be the coarse mesh computation, as was typically
done in the multifidelity method {\color{blue}\cite{ZNX14, ZX17}}.

This paper is organized as follows. Section \ref{sec:1} gives a general introduction on the bi-fidelity stochastic collocation method for PDEs with random parameters. 
{\color{blue} In Section \ref{sec:2A}, we introduce the Boltzmann equation and its fluid limit, and the corresponding high- and low-fidelity models. In section \ref{sec:2} we establish a {\it uniform-in-$\e$} error estimates of this bi-fidelity modeling for the Boltzmann equation. We conduct the error analysis in two cases, one for all regimes and one particularly for the fluid regimes (where $\e \ll 1$). In Section \ref{sec:3}, other possible choices of low-fidelity models are discussed, in particular, the 
low fidelity model being the same equation but solved by using a coarser mesh.} As examples, we show the error estimates of discretized asymptotic-preserving schemes 
for the uncertain linear transport equations in the diffusive limit. The paper is concluded in Section \ref{sec:4}.

\section{A bi-fidelity stochastic collocation method}
\label{sec:1}

We first briefly review the efficient bi-fidelity stochastic collocation method studied in \cite{NGX14, ZNX14}. {\color{red} We denote the high- and low-fidelity solutions by $u^H$, $u^L$ respectively, while $u^B$ is the sought bi-fidelity solution, which is an approximation of $u^H$. }
The bi-fidelity algorithm consists of two stages. In the offline stage, one employs the cheap low-fidelity model to explore the random parameter space. One selects the most important parameter points $\gamma_N$ by the greedy procedure as a popular reduced-basis method for solving the parameterized partial differential equations \cite{RB89}. 

Within the online stage, the bi-fidelity approximation is realized by applying exactly the same approximation rule learnt from the low-fidelity model for any given $z$.  For any given sample point $z\in I_z$, we project the low-fidelity solution $u^L(z)$ onto the low-fidelity approximation space
$U^L(\gamma_N)$: 
\begin{equation}\label{PU}  u^L(z) \approx \mathcal P_{U^L(\gamma_N)}[u^L(z)] = \sum_{k=1}^N c_k(z) u^L(z_k), \end{equation}
where $\mathcal P_{U^L(\gamma_N)}$ is the projection operator onto the space $U^L(\gamma_N)$ with the corresponding projection coefficients $\{ c_k\}$ computed by the Galerkin approach
\begin{equation}\label{Gc}
{\bf G}^L {\bf c} = {\bf f}, \qquad {\bf f} = (f_k)_{1\leq k\leq N}, \qquad f_k = 
\langle u^L(z), u^L(z_{k})\rangle^L. 
 \end{equation}
Here ${\bf G}^L$ is the Gramian matrix of $u^L(\gamma_N)$, 
\begin{equation}\label{GM}
 ({\bf G}^L)_{ij} =  \left\langle u^L(z_i), u^L(z_j) \right\rangle^L, \qquad 1 \leq i,\, j \leq N, 
\end{equation}
with $\langle\cdot,\cdot\rangle^L$ the inner product associated with the approximation space $U^L(\gamma_N)$.
These low-fidelity coefficients $\{c_k\}$ serve as the surrogate of the corresponding high-fidelity coefficients of $u^H(z)$. 
The bi-fidelity approximation of $u^H$ can be constructed as follows: 
\begin{equation}\label{UB}
 u^B(z) = \sum_{k=1}^N c_k(z) u^H(z_k).  
\end{equation}
The key idea of this method is that one first writes the solution $u^L(z)$ as coordinates in the basis $u^L(\gamma)$, which is assumed to be a collection 
of linearly independent solutions, followed by using exactly the same coordinates coefficients $c_n$ in the reconstruction for the bi-fidelity solution \eqref{UB}. 
In practice, the number of low-fidelity basis is typically small, the cost of computing the projection coefficients ${\bf c}$ by solving the linear system \eqref{Gc} is negligible, thus the dominant cost of the online step is only one low-fidelity simulation run. We emphasize that the low-fidelity coefficients would be a good approximation of the corresponding high-fidelity coefficients for a given sample $z$, 
if the low-fidelity model can mimic the variations of the high-fidelity model in the parameter space. 

Details of this bi-fidelity approximation are summarized in Algorithm \ref{alg} below. Though most of the steps are straightforward, we refer readers to \cite{hampton2018practical, NGX14, ZNX14} for more technical details and justifications. 
\\[3pt]
\begin{algorithm}[H]
\caption{bi-fidelity approximation}
\label{alg}
\label{BiFi-pod}
\Offline

Select a sample set $\Gamma = \{z_1, z_2, \hdots, z_M\}\subset I_z $.

Run the low-fidelity model $u_l(z_j)$ for each $z_j \in \Gamma$.

Select $N$ ``important" points from $\Gamma$ and denote it by $\gamma_N=\{z_{i_1}, \cdots z_{i_N} \} \subset\Gamma$. 
Construct the low-fidelity approximation space $U^L(\gamma_N)$.

Run high-fidelity simulations at each sample point of the selected sample set $\gamma_N$. Construct the high-fidelity approximation space $U^H(\gamma_N)$. 

\Online

For any given $z$, get the low-fidelity solution $u^L(z)$ and compute the low-fidelity coefficients by projection:
\begin{equation}\label{C-N}
u^L(z) \approx \mathcal{P}_{U^L(\gamma_N)}u^L = \sum_{k=1}^N c_k(z)u^L(z_k). \end{equation}

Construct the bi-fidelity approximation by applying the sample approximation rule in the low-fidelity model:
$$u^B(z)  = \sum_{k=1}^N c_k(z)u^H(z_k).$$
\end{algorithm}

In order to obtain the error estimate of $u^H - u^B$ in general, we use the following way to split the total error, by inserting the information 
of $u^L$: 
\begin{align}
\label{idea}
\begin{split}
&\displaystyle\quad u^H(z) - u^B(z)  \\[4pt]
&\displaystyle = u^H(z) - \sum_{n=1}^N c_n(z) u^H(z_n)  \\[4pt]
&\displaystyle = u^H(z) - u^L(z) + \left(u^L(z) -  \sum_{n=1}^N c_n(z) u^L(z_n)\right) + \sum_{n=1}^N c_n(z)\left( u^L(z_n) - u^H(z_n)\right), 
\end{split}
\end{align}
where the second term is nothing but the projection error of the greedy algorithm, and it remains to estimate $u^H(z) - u^L(z)$ in proper norms. 

\section{Bi-fidelity modeling for the Boltzmann equation}
\label{sec:2A}

In this section, we introduce a bi-fidelity modeling for the Boltzmann equation using the framework in Section \ref{sec:1}. 

There are usually small parameters in the kinetic equations, characterized by the Knudsen number (mean free path), defined as the average distance between two collisions of particles. When the Knudsen number (denoted by $\e$) goes to $0$, the kinetic equations pass to the limit of their hydrodynamic equations, which describe the macroscopic behavior and own cheaper computational cost, thus can be chosen as the efficient low-fidelity models. 
However, when $\e=O(1)$, although the macroscopic models may no longer be close to the kinetic ones, thus the error estimates cannot be obtained easily, we can still obtain an error estimate in the {\it perturbative setting}--in which the solutions are close to the global equilibrium, using the hypocoercivity argument 
\cite{MB15,LJ-UQ}. 

We first give an introduction to the Boltzmann equation,  the most celebrated kinetic equations, for rarefied gas. 
There are many sources of uncertainties in the Boltzmann equation, such as
the initial data, boundary data, and collision kernel. A dimensionless form with scalings reads
\begin{align}
\label{Boltz}
\left\{
\begin{array}{l}
\displaystyle\partial_{t} f+\frac{1}{\varepsilon^{\alpha}} v \cdot \nabla_{x} f=\frac{1}{\varepsilon^{1+\alpha}} \mathcal{Q}(f, f),  \\[8pt]
\displaystyle  f(0, x, v, z)=f_{\text{in}}(x, v, z), \quad x \in \Omega \subset \mathbb{T}^{d_x}, v \in \mathbb{R}^{d_v}, z \in I_{z}  \subset \mathbb R^{d_z}, 
\end{array}\right.
\end{align}
where $f(t,x,v,z)$ is the probability density distribution function, modeling
the probability of finding a particle at time $t$, position $x\in \Omega \subset \mathbb{T}^{d_x}$ (periodic box of $d_x$ dimension), and with velocity 
$v \in \mathbb{R}^{d_v}$. We assume the volume of the spatial domain $\Omega$ is bounded. 
The parameter $\e$ is the Knudsen number. 
The case $\alpha=0$ corresponds to the acoustic scaling, and $\alpha=1$ corresponds to the 
incompressible Navier-Stokes (INS) scaling. 
The collision operator $\mathcal Q$ is a quadratic integral operator modeling the binary elastic collision between particles, and is given by 
\begin{equation} \mathcal Q(f, f)(v) = \int_{\mathbb R^{d_v}}\int_{\mathbb S^{d_v -1}} 
B(|v-v_{\ast}|, \cos\theta, z) \left(f(v^{\prime})f(v_{\ast}^{\prime}) - f(v)f(v_{\ast})\right) d\sigma dv_{\ast}. 
\end{equation}
The velocity pairs before and after the collision $(v, v_{\ast})$ and $(v^{\prime}, v_{\ast}^{\prime})$ have the relation: 
\begin{align}
\begin{cases}
&\displaystyle v^{\prime}= \frac{v+v_{\ast}}{2} + \frac{|v-v_{\ast}|}{2}\sigma, \\[4pt]
&\displaystyle v_{\ast}^{\prime} = \frac{v+v_{\ast}}{2} - \frac{|v-v_{\ast}|}{2}\sigma, 
\end{cases}
\end{align}
with the vector $\sigma$ the scattering direction varying on the unit sphere $\mathbb S^{d-1}$. 
The collision kernel $B$ is a non-negative function depending on the modulus of
the relative velocity $|v-v_{\ast}|$, cosinus of the deviation angle $\theta$ with $$\cos\theta = \sigma \cdot (v-v_{\ast})/|v-v_{\ast}|, $$ and the random variable 
$z\in I_z\subset \mathbb R^{d_z}$. Periodic boundary condition is considered. 

Denote $m(v) = \left(1, v, \frac{|v|^2}{2}\right)^{\mathbb T}$, then 
\begin{equation}\label{Q_cons}\int_{\mathbb R^{d_v}}\mathcal Q(f, f)m(v)\, dv = 0, \end{equation}
which correspond to the conservations of mass, momentum and total energy of the collision operator. 
The celebrated Boltzmann's H-theorem \cite{Cercignani-Book} gives: 
$$ \int_{\mathbb R^d} \mathcal Q(f, f)\ln f\, dv \leq 0. $$
The equality holds if and only if $f$ reaches the equilibrium state known as the local Maxwellian: 
\begin{equation}\label{Maxwell} M(v)_{\rho, u, T} = \frac{\rho}{(2\pi T)^{\frac{d_v}{2}}}\exp\left(-\frac{|v-u|^2}{2T}\right), 
\end{equation}
where $\rho$, $u$ and $T$ are the density, bulk velocity and temperature, respectively: 
\begin{equation}\label{macro}\rho = \int_{\mathbb R^{d_v}} f(v)\, dv, \qquad u = \frac{1}{\rho}\int_{\mathbb R^{d_v}} f(v)v\, dv, 
\qquad T = \frac{1}{d_v \rho} \int_{\mathbb R^{d_v}} f(v)|v-u|^2\, dv. \end{equation}

We consider the case of hard potential and Maxwellian molecules, that is, the collision kernel takes the form
 \begin{align}
 \label{B-Form} 
 &\displaystyle B(|v-v^{\ast}|, \cos\theta, z) = \Phi(|v-v^{\ast}|)\,b(\cos\theta, z), \qquad b(\cos\theta, z) >0,  \\[4pt]
 &\displaystyle\label{Phi}\Phi(|v-v^{\ast}|) = C_{\Phi} |v-v^{\ast}|^{\gamma},  \qquad C_{\Phi}>0,\, \gamma \in [0,1]. 
 \end{align}
 \hspace{2cm}
 
\noindent{\bf \color{blue}The fluid limit}\, 
We introduce the fluid approximation to the Boltzmann equation when $\e \to 0$, known as the 
compressible Euler system. Denote $\langle\, \cdot\, \rangle$ as the velocity averages of the argument, 
$$ \langle f \rangle = \int_{\mathbb R^{d_v}} f(v)\, dv. $$
Multiplying (\ref{Boltz}) by $m(v)$ and integrating with respect to $v$, by the conservation property of $\mathcal Q$ given in \eqref{Q_cons}, 
and the approximation $f\to M(v)_{\rho, u, T}$ when $\e \to 0$, one gets the compressible Euler equations
of gas dynamics:
\begin{equation}\label{Euler}
\partial_t \begin{pmatrix} \rho \\ \rho u \\ E \end{pmatrix} + 
\nabla_x \cdot \begin{pmatrix} \rho u \\ \rho u\otimes u + p\,\text{I} \\ (E+p) u \end{pmatrix} = 0, 
\end{equation}
where $E$ and $T$ given in \eqref{macro} have the relation: 
\begin{equation}\label{E-T}
E= \frac{1}{2} \rho\, |u|^2 +\frac{d_v}{2} \rho T. \end{equation}

{\color{red} Since our goal is to efficiently solve the Boltzmann equation, we always let the uncertain Boltzmann equation \eqref{Boltz} be the high-fidelity model.}
As implemented in \cite{LZ19}, we now choose the compressible Euler system \eqref{Euler} as our low-fidelity model.
It is a first-order $O(\e)$ approximation to the Boltzmann equation, and can mimic the variations in the random space of macroscopic quantities of the Boltzmann equation up to a certain accuracy, see numerical examples in \cite{LZ19}. 

Let $u^H$ be the high-fidelity solution, which is defined by the $d_v+2$ macroscopic moments of density, momentum and energy that are obtained from the distribution $f$ solved by the Boltzmann equation: 
$$ u^H = \int_{\mathbb{R}^{d_v}}\left(\begin{array}{c}{1} \\ {v} \\ {\frac{1}{2}|v|^{2}}\end{array}\right) f(v) dv := \left(\begin{array}{c}{\rho^H} \\ \rho^H u^H \\ E^H\end{array}\right) := \left(\begin{array}{c} u_1^H \\ u_2^H \\ u_3^H \end{array}\right). $$

\section{Error estimate of the bi-fidelity modeling for the Boltzmann equation}
\label{sec:2}

In this section we conduct the error analysis for the bi-fidelity method introduced in previous section. The analysis is based on the hypocoercivity analysis for uncertain Boltzmann equation, as introduced in \cite{LJ-UQ}.

\subsection{The mathematical setting}

{\color{red}In order to adopt the hypocoercivity theory and study long-time behavior of the solution to \eqref{Boltz}, 
one has to consider the perturbative setting. In this framework, we let $h$ be the perturbed solution, then consider the ansatz} \begin{equation}\label{PS} f = M + \delta \sqrt{M}\, h, \end{equation}
{\color{red} where parameter $\delta$ is assumed sufficiently small and independent of $\varepsilon$.}
Denote the steady state 
by $\displaystyle u^{st}=(1,0,\frac 3 2)^T = (u_1^{st}, u_2^{st}, u_3^{st})$, then 
\begin{equation}\label{uH-def}  u^H = u^{st} + \delta \int_{\mathbb R^{d_v}}
\left(\begin{array}{c}{1} \\ {v} \\ {\frac{1}{2}|v|^{2}}\end{array}\right) \sqrt{M}h(v)  dv. \end{equation}

We let the low-fidelity solution 
\begin{equation} u^L = \left(\begin{array}{c}\rho^L, \rho^L u^L, E^L\end{array}\right) := \left(\begin{array}{c} u_1^L \\ u_2^L \\ u_3^L \end{array}\right), 
\end{equation} obtained from the compressible Euler system \eqref{Euler}. 
Since the initial condition of the high-fidelity and low-fidelity models are required to be consistent,  
$$\label{uL-IC} u^L(t=0) = \int_{\mathbb{R}^{d_v}}\left(\begin{array}{c}{1} \\ {v} \\ {\frac{1}{2}|v|^{2}}\end{array}\right) f_{\text{in}}(v) dv. $$
We then consider a linearization of the Euler equations \eqref{Euler} around the state $(\rho_0, u_0, T_0) = (1,0,1)$, 
\begin{equation}\label{Rho-E} \rho^L = 1+ \delta\, \tilde\rho, \qquad u^L = \delta\, \tilde u, \qquad T^L = 1 + \delta\, \tilde T, \end{equation}
{\color{red} where $\delta$ is the same  as that in the high-fidelity solution expansion \eqref{uH-def}, since we want the initial data of $u^L$ and $u^H$ to be consistent.}
One can derive the following acoustic equations: 
\begin{equation}
\label{AC-eqn}
\left\{
\begin{array}{lll}
\displaystyle\partial_{t} \tilde\rho+\nabla_{x} \cdot \tilde u = 0, \\[4pt]
\displaystyle\partial_{t}\tilde u+\nabla_{x}(\tilde\rho+\tilde T) = 0, \\[4pt]
\displaystyle\frac{d_v}{2}\, \partial_{t}\tilde T +\nabla_{x} \cdot\tilde u = 0, 
\end{array}
\right. 
\end{equation}
which is essentially a wave equation \cite{Bardos2000,Russel}, and obtaining $E^L$ from $T^L$ is shown in \eqref{E-T}. 

We first define the space and norms that will be used. 
The Hilbert space of the random variable is given by 
$$H\left(\mathbb{R}^{d} ; \pi\,\mathrm{d} z\right)=\left\{f | I_z \rightarrow \mathbb{R}, \int_{I_z} f^{2}(z) \pi(z) \mathrm{d} z< \infty\right\}, $$
and equipped with the inner product  
$$\langle f, g\rangle_{\pi}=\int_{I_z} f g \,\pi(z) \mathrm{d}z. $$ 

{\color{red} Define the Sobolev norm of $h$ by }
$$ || h ||_{H_{x,v}^s}^2 := \sum_{|m|+|n|\leq s} ||\partial_x^m \partial_v^n h ||_{L^2_{x,v}}^2, $$
and the summation of Sobolev norms for $z$-derivatives of $h$, 
$$\|h\|_{H_{x, v}^{s, r}}^{2}=\sum_{|m| \leq r}\left\|\partial^{m} h\right\|_{H_{x, v}^{s}}^{2}, $$
in addition to the norms in the $(x,v,z)$--space:
$$\|h\|_{H_{x, v}^{s} H_{z}^{r}}^{2}=\int_{I_{z}}\|h\|_{H_{x, v}^{s,r}}^2\, \pi(z) d z, \qquad
 \|h\|_{H_{x, v}^{s} L_z^2}^{2}=\int_{I_{z}}\|h\|_{H_{x, v}^{s} L_z^2}^2\, \pi(z) d z, $$
and the sup norm in $z$: 
$$\|h\|_{H_{x, v}^{s, r} L_{z}^{\infty}}=\sup _{z \in I_{z}}\|h\|_{H_{x, v}^{s, r}}. $$
We introduce the standard multivariate notation. Denote the countable set of
``finitely supported" sequences of nonnegative integers by
$$
\mathcal{F} :=\left\{\nu=\left(\nu_{1}, \nu_{2}, \cdots\right) : \nu_{j} \in \mathbb{N}, \text { and } \nu_{j} \neq 0 \text { for only a finite number of } j\right\}
$$
with \(|\nu| :=\sum_{j \geq 1}\left|\nu_{j}\right| .\) For \(\nu \in \mathcal{F}\) supported in \(\{1, \cdots, J\},\) the partial
derivative in \(z\) is defined by 
$$
\partial_z^{\nu} u=\frac{\partial^{|\nu|} u}{\partial^{\nu_{1}} z_{1} \cdots \partial^{\nu_{J}} z_{J}}\,.
$$

We make the following assumptions on the random collision kernel and initial data:  
\begin{assumption}\label{Assump}
Assume that each component of the random variable $z :=\left( z_j \right)_{j \geq 1}$ has a compact support. 
The collision kernel satisfies 
\begin{equation} \label{Coll} 0 < b(\mu, z) \leq C_b, \quad |\partial_{\mu}b(\mu, z)| \leq \widetilde C_b,  \quad |\partial_z^{\nu} b(\mu, z)| \leq C, \end{equation}
where $\mu=\cos\theta \in [-1,1]$, and $|\nu| \leq r$ (the constant $r$ is associated to the regularity of the initial data in the random space), $C_b$, $\widetilde C_b$ and $C$ are all constants. 
Let $\left(\psi_j \right)_{j\geq 1}$ be an {\it affine representer} of the random initial data $h_{\text{in}}$, 
which by definition means that \cite{Cohen15}
 \begin{equation}\label{h-AF} h_{\text{in}}(z) = \tilde h_0 + \sum_{j \geq 1} z_j \psi_j, 
 \end{equation}
 where $\tilde h_0 = \tilde h_0(x,v)$ is independent of $z$, and the sequence $\left( ||\psi_j||_{L^{\infty}(V)}\right)_{j\geq 1}\in \ell^p$ for $0<p<1$, with $V$ representing the physical space. 
 \end{assumption}
{\color{red} The reason we need to assume compact support of random variable $z$ and \eqref{h-AF} is because one needs such an assumption in  \cite{Cohen15}}, the result we rely upon when we estimate the projection 
error for greedy algorithm later this section. This is only a technical assumption while the numerical implemendation does not require this.

\subsection{The analysis for all $\e$}

{\color{red}We first show the main result of this subsection: }
\begin{theorem}
\label{Thm1}
Let \textbf{Assumption 1} hold and assume the initial data satisfies 
\begin{equation}\label{h-IC} || h_{\text{in}} ||_{H_{x,v}^{s,r} L_z^{\infty}} \leq \frac{\varepsilon}{\delta}\,\eta, \end{equation}
for sufficiently small $\eta$. 
By implementing $N$ high-fidelity simulation runs, the error estimate is given by
\begin{equation}\label{EE} \left\| u_i^H(z)  - u_i^B(z) \right\|_{H_{x}^{s}L_z^2} \leq \frac{C_1}{(N/2+1)^{q/2}} + \frac{1}{\sqrt{\lambda_0}} \max\{\eta^{\prime}, C_{\xi}\} \left(1 + e^{-\e^{1-\alpha}\tau t}\right)\sqrt{N}\delta, \end{equation}
where $C_1$, $C_2$ are independent of $\e$ and $N$, and depend on the initial data of the perturbation solutions $\eta$ and $\xi$, which are assumed 
sufficiently small. 
\end{theorem}
\hspace{1cm}

\noindent {\textbf{Proof.}} Under the assumption for the collision kernel \eqref{Coll}, with some modification to \cite{LJ-UQ}, it can be derived that if the initial data satisfies \eqref{h-IC}, 
then at all time $t>0$, 
\begin{equation}\label{h-eps}\left\|h_{\varepsilon}\right\|_{H_{x, v}^{s, r} L_{z}^{\infty}} \leq \eta^{\prime}\, e^{-\e^{1-\alpha}\tau t}, \quad \text{ and   }\left\|h_{\varepsilon}\right\|_{H_{x, v}^{s} H_{z}^{r}} \leq \eta^{\prime}\, e^{-\e^{1-\alpha}\tau_{s} t}. \end{equation}
To not distract the reader, we show its proof in Appendix A. We would like to remark that we let $\delta$ in the perturbative setting \eqref{PS} to be independent of $\e$, while 
$\delta=\e$ is enforced in the deterministic analysis \cite{MB15}. {\color{black} Combined \eqref{h-IC} with \eqref{h-AF} in {\bf Assumption 1}, we actually only need $r=1$.}
We carry out the following estimate componentwise for $u^H$ and $u^L$. 
By \eqref{uH-def} and using the Cauchy-Schwarz inequality for each $i$, one gets 
\begin{equation}\label{u-H} ||u_i^H - u_i^{st} ||_{H_{x}^{s}H_z^r} \leq \delta\, ||h||_{H_{x}^{s}L_v^2 H_z^r} \leq 
\delta\, ||h||_{H_{x,v}^s H_z^r} \leq \delta\, \eta^{\prime}\, e^{ - \e^{1-\alpha}\tau_s t}, \quad i=1,2,3. \end{equation}
\hspace{2cm}

Since the acoustic equations \eqref{AC-eqn} is a system of linear hyperbolic equations, it is clear by the method of characteristics that 
if $||\tilde u_i^L(t=0)||_{H_{x}^{s}H_z^r} \leq \xi$ then at all time $t>0$, 
\begin{equation}\label{E-Rho} ||\tilde u^L_i||_{H_{x}^{s}H_z^r} \leq C_{\xi}, \end{equation}
where $\tilde u^L_i$ denotes the perturbed part of the low-fidelity solution $u^L_i$. 
Due to \eqref{uL-IC} and \eqref{h-IC}, the bound $\xi$ is also sufficiently small. 
From \eqref{E-Rho}, we now have 
 \begin{equation}\label{u-L0}  || u^L_i - u^{st}_i||_{H_{x}^{s} H_z^r} = \delta\, ||\tilde u^L_i||_{H_{x}^{s} H_z^r} \leq C_{\xi}\, \delta. \end{equation}
If the initial data satisfies
\begin{equation}\label{H-IC} || h_{\text{in}} ||_{H_{x,v}^{s,r}L_z^{\infty}}\leq \frac{\varepsilon}{\delta}\,\eta, 
 \end{equation}
 then by \eqref{u-H} and \eqref{u-L0}, 
\begin{align}
\label{u-HL0}
\begin{split}
||u^H_i - u^L_i||_{H_{x}^{s}H_z^r} & = ||u^H_i-u^{st}_i -(u^L_i-u^{st}_i)||_{H_{x}^{s} H_z^r}  \\[4pt]
& \leq  ||u^H_i - u^{st}_i ||_{H_{x}^{s}H_z^r}  + || u^L_i - u^{st}_i||_{H_{x}^{s} H_z^r} \\[4pt]
& \leq \eta^{\prime}\,\delta\, e^{ -\e^{1-\alpha}\tau_s t} + C_{\xi}\,\delta, 
\end{split}
\end{align}
therefore, 
\begin{equation}\label{u-HL}
||u^H_i - u^L_i||_{H_{x}^{s}L_z^2} \leq  \max\{\eta^{\prime},\, C_{\xi}\}\, \delta\, (1 + e^{ -\e^{1-\alpha}\tau_s t}) \leq C^{\prime}\, \delta\,. 
\end{equation}
This indicates that the difference of density, momentum and energy obtained between the high- and low-fidelity solutions are all of $O(\delta)$. 
\begin{remark}
\label{Rmk}
We make a remark that the estimate \eqref{u-HL} still holds in the kinetic regime when $\e = O(1)$. Our error estimate \eqref{u-HL0} is
{\it uniform in $\e$}. We indeed require that $\delta \ll 1$, i.e., 
the classical solution to the Boltzmann equation \eqref{Boltz} is assumed to be near the global equilibrium, with sufficiently small initial data 
of the perturbed solution. We use both the regularity of low-fidelity and high-fidelity solutions in the proof. 
Actually, in order to get the result \eqref{u-HL} for the Boltzmann equation, it is not necessary to choose hydrodynamic limit models, compressible Euler system in our case, as the low-fidelity models. We only require the perturbative settings \eqref{PS} and \eqref{Rho-E} be consistent, meaning that the $d_v+2$ velocity moments of $f$ (or $u^H$) shown in \eqref{PS} satisfies exactly the same perturbative setting as that for the low-fidelity solution $u^L$ shown in \eqref{Rho-E}, namely, the moments of the high-fidelity kinetic model are $O(\delta)$ away from the low-fidelity variables. This suggests that any moment closure model, as long as it has moments approximately closed to the Boltzmann equation, can serve as the low-fidelity model and our analysis applies. 
\end{remark}

Take $||\cdot||_{H_x^s L_z^2}$ norm on both sides of the equality \eqref{idea}, for each moment component $i=1,2,3$, we get
\begin{align}
\label{idea-1}
\begin{split}
&\quad \left\| u^H_i(z) - u^B_i(z) \right\|_{H_{x}^{s}L_z^2} =\left \|u^H_i(z) - \sum_{n=1}^N c_n(z)\, u^H_i(z_n)\right\|_{H_{x}^{s}L_z^2}  \\[4pt]
& \leq \left\| u^H_i(z) - u^L_i(z) \right\|_{H_{x}^{s}L_z^2} +   \left\| u^L_i(z) -  \sum_{n=1}^N c_n(z)\, u^L_i(z_n)\right\|_{H_{x}^{s}L_z^2} 
+ \underbrace{\left\|\, \sum_{n=1}^N c_n(z) \left(u^L_i(z_n) - u^H_i(z_n) \right) \right\|_{H_{x}^{s}L_z^2}}_{\text{Term A}}. \end{split}
\end{align}
The first term is $O(\delta)$ by \eqref{u-HL}, and the second term is actually the projection error of the greedy algorithm when searching the most important points $\gamma_N$ from the low-fidelity solution manifold. To estimate the third term, we first get a bound for the vector $||{\bf c}||$, with $||\cdot||$ the matrix induced $\ell_2$ norm. 

From the definition of projection onto $U_i^L(\gamma_N)$ for each $u_i^L(z)$ with $i=1,2,3$, $\mathcal P_{U_i^L(\gamma_N)}[u_i^L(z)]$ given in \eqref{PU}, 
$$ \left(\mathcal P_{U^L(\gamma_N)}[u_i^L(z)] \right)^2 = \sum_{m,n=1}^N c_m(z)c_n(z)u_i^L(z_m)u_i^L(z_n), $$
thus \begin{align*}\displaystyle \int_{\Omega}\left(\mathcal P_{U_i^L(\gamma_N)}[u_i^L(z)] \right)^2 dx 
& = \sum_{m,n=1}^N c_m(z)c_n(z)\int_{\Omega}u_i^L(z_m)u_i^L(z_n) dx \\[4pt]
&\displaystyle  := {\bf c}^T {\bf G}^{L} {\bf c} \geq \lambda_0 \left\|{\bf c}\right\|^2, 
\end{align*}
where ${\bf G}^{L}$ is the Gramian matrix of $u_i^L(\gamma_N)$ defined in \eqref{GM} (where $\left\langle \cdot, \cdot \right\rangle$ applies to $L_x^2$ here), and $\lambda_0>0$ is its minimum eigenvalue (we omit the $i$ in ${\bf G}^L$, $\lambda_0$ and ${\bf c}$ for $i=1,2,3$ here). 
The last inequality is due to Gramian matrices are positive-semidefinite and we assume it positive-definite (otherwise ${\bf c}$ could not be obtained by solving the system \eqref{Gc}). 
Since for all $z$,  $$\int_{\Omega}\left(\mathcal P_{U^L(\gamma_N)}[u^L(z)] \right)^2 dx \leq \int_{\Omega} [u^L(z)]^2\, dx, $$  then 
\begin{equation}\label{CC}
 ||{\bf c}|| \leq \frac{1}{\sqrt{\lambda_0}} \left(\int_{\Omega}\left(\mathcal P_{U_i^L(\gamma_N)}[u_i^L(z)] \right)^2\, dx \right)^{1/2} \leq 
\frac{1}{\sqrt{\lambda_0}} ||u_i^L(z)||_{L^2_x}, 
\end{equation}
for all $z$. Thus $\displaystyle\Big\| ||{\bf c}|| \Big\|_{L_z^{\infty}} \leq \frac{1}{\sqrt{\lambda_0}} ||u_i^L(z)||_{L^2_x L_z^{\infty}}$. 
Since $\sqrt{x}$ is a monotone function, $\Big\| ||{\bf c}|| \Big\|_{L_z^{\infty}}=\left( \sum_{n=1}^N ||c_n(z)||_{L_z^{\infty}}^2 \right)^{1/2}$. 
By the regularity of $u_i^L$ in \eqref{u-L0}, and the assumption that the volume of $\Omega$ is bounded, \eqref{CC} implies that
\begin{equation}\label{C-N} \left( \sum_{n=1}^N ||c_n(z)||_{L_z^{\infty}}^2 \right)^{1/2} \lesssim  \frac{1}{\sqrt{\lambda_0}} || u_i^L(z) ||_{L_x^2 L_z^{\infty}}
\leq \frac{1}{\sqrt{\lambda_0}} || u_i^L(z) ||_{H_x^s L_z^{\infty}} 
\lesssim \frac{1}{\sqrt{\lambda_0}}(u_i^{st} +C_{\xi}\,\delta). \end{equation}

By the Cauchy-Schwarz inequality, using \eqref{u-HL} and \eqref{C-N}, one gets 
\begin{align}
\label{T-A}
\begin{split}
\text{Term A} &\leq \left( \sum_{n=1}^N ||c_n(z)||_{L_z^2}^2 \right)^{1/2}  \left( \sum_{n=1}^N ||u_i^H(z_n) - u_i^L(z_n)||_{H_x^s L_z^2}^2\right)^{1/2} 
\\[4pt]
&\lesssim \left( \sum_{n=1}^N ||c_n(z)||_{L_z^{\infty}}^2 \right)^{1/2} \left( \sum_{n=1}^N ||u_i^H(z_n) - u_i^L(z_n)||_{H_x^s}^2\right)^{1/2} \\[4pt]
& \leq \sqrt{N}  \left( \sum_{n=1}^N ||c_n(z)||_{L_z^{\infty}}^2 \right)^{1/2} \left( \max_{n}||u_i^H(z_n) - u_i^L(z_n)||_{H_x^s} \right) \\[4pt]
& \leq \frac{C^{\prime}}{\sqrt{\lambda_0}} \sqrt{N}\, \delta\,  (u_i^{st} + C_{\xi}\,\delta) \lesssim C \sqrt{N}\, \delta, 
\end{split}
\end{align}
where $C^{\prime}$ and $C$ are generic constants. The boundness of the random variable $z$ is used in the first inequality. 

Now we estimate the second term on the right-hand-side of \eqref{idea-1} below (and omit the subscript $i$ for $i=1,2,3$ in $u^L(z)$). 
The Kolmogorov $N$-width of a functional manifold describes the ``best" achievable distance for approximation from a general $N$-dimensional subspace. We denote $d_N(u^L(I_{z}))$ the Kolmogorov $N$-width of the functional manifold $u^L(I_{z})$, defined by
$$d_N (u^L(I_{z})) = \inf_{\text{dim}(V_N)=N} \sup_{v\in u^L(I_{z})} d^L(v, V_N). $$
 Denote the space $\mathcal H=H_{x}^s$. From for example \cite{Cohen15,Devore13}, we recall that the projection error of the greedy algorithm satisfies
\begin{equation}\label{NW} 
\sup_{z\in I_{z}} \left\| u^L(z) - P_{U^L(\gamma_N)}u^L(z)\right\|_{\mathcal H} \leq C \sqrt{d_{N/2}(u^L(I_{z}))}\,. 
\end{equation}
In \cite[Section 4.2]{LZ19}, the upper bound for the Kolmogorov $N$-width is obtained based on \cite{Cohen15}.
Let the $N$-dimensional subspace $V_N := \text{Span}\left\{ c_{\nu} : \nu\in\Lambda_N\right\}$. 

The analyticity of $u^L$ is obvious.
From \eqref{h-AF}, due to the affine dependence on $z$ of the initial data $h_{\text{in}}$,       
it is analytic in the random space.  
Since \eqref{AC-eqn} is a {\it linear} acoustic system, thus its solution at all time stays analytic with respect to the randomness. 
Similar to \cite[Section 4.2]{LZ19}, one has
\begin{align}
\label{d-N}
\begin{split} d_{N}(u^L(I_z))_{\mathcal H} & \leq \sup _{v \in u^L(I_z)} \min _{w \in V_{N}}\|v-w\|_{\mathcal H} =\sup _{z \in I_{z}} \min _{w \in V_{N}}\left\|u^L(z)-w\right\|_{\mathcal H} \\ & \leq\left\|u^L-\sum_{\nu \in \Lambda_{N}} w_{\nu} P_{\nu}\right\|_{L^{\infty}\left(I_{z}, \mathcal H\right)} \leq C(N+1)^{-q}, \quad q=\frac{1}{p}-1. 
\end{split}
\end{align}
The constant $p$ is associated to the affine representer $\left(\psi_j \right)_{j\geq 1}$ in the random initial data shown in \eqref{h-AF}. 
The sum $\sum_{\nu\in\Lambda_N} w_{\nu}P_{\nu}$ is the truncated Legendre expansion with $\left(P_k\right)_{k \geq 0}$ the sequence of
renormalized Legendre polynomials on $[-1,1]$, and $\Lambda_N$ is the set of indices that corresponds to the $N$ largest $||w_{\nu}||_{\mathcal H}$.  
\begin{remark}
Note that compared to \cite[Theorem 4.1]{LZ19} where 
$$\sup_{z\in I_{z}} ||u^H(z) - P_{U^H(\gamma_N^L)}u^H(z)|| \leq C \sqrt{d_{N/2}(u^H(I_{z}))}$$
 is employed in their estimate, it is different here since we directly use the projection error of the greedy algorithm exploring from the low-fidelity 
 solution manifold $u^L(z)$. 
 Therefore we don't need to prove the analyticity of the high-fidelity solution solved by the Boltzmann equation, which seems a quite challenging task. However, for the sake of future research, we prove an analytic result in Appendix B for the linearized Boltzmann equation.
\end{remark}

By \eqref{idea-1}, \eqref{T-A}, \eqref{NW}, \eqref{d-N} and note that $N/2$ needs to be plugged into \eqref{d-N} due to \eqref{NW}, we obtain \eqref{EE} and proved Theorem \ref{Thm1}. 
\qquad\qquad\qquad\qquad\qquad $\square$

\begin{remark}
In typical numerical experiments carried out in \cite{LZ19}, $N\approx 20$. If additionally one assumes that $\delta = o(1/\sqrt{N})$, 
then \begin{equation}\label{Result1} \left\| u_i^H(z)  - u_i^B(z) \right\|_{H_{x}^{s}L_z^2} \leq \frac{C_1}{(N/2+1)^{q/2}} + C_2, \end{equation}
where $C_2=o(1)$, and is independent of $N$ and $\e$. 
\end{remark}
\subsection{Analysis for $\e \ll 1$}
\label{analy2}

When $\e \ll 1$, we incorporate the convergence of the perturbative solution $h_{\e}$ towards its limit \cite[Section 8]{MB15}, then get a sharper error estimate of the bi-fidelity method. 
Under the acoustic scaling $\alpha=0$ in \eqref{Boltz}, $(\rho, u, \theta)$ satisfy the acoustic system \eqref{AC-eqn}, {\color{red} which is chosen as the low-fidelity model $u^L$.}
Under the INS scaling $\alpha=1$ in \eqref{Boltz}, they are weak solutions \cite{Leray}
to the incompressible Navier-Stokes equations, {\color{red}which is considered as the low-fidelity model $u^L$ in this case: }
 \begin{equation}
\left\{
\begin{array}{lll}
\displaystyle\partial_{t} u-\nu \Delta u+u \cdot \nabla u+\nabla p=0, \\[2pt]
\displaystyle\nabla \cdot u=0, \\[2pt]
\displaystyle\partial_{t}\theta -\kappa \Delta\theta +u \cdot\nabla\theta =0, 
\end{array}
\right. 
\end{equation}
where $p$ is the pressure, $\nu$ and 
$\kappa$ are constants, together with the Boussinesq relation $\nabla(\rho + \theta)=0$. 

The main result of this section is summarized below: 
\begin{theorem}
\label{Thm2}
If the initial data $h_{\text{in}}$ belongs to $H_x^{s+\ell}L_v^2 L_z^{\infty}$ where $\ell\in (0,1]$, then at any time $T>0$, 
\begin{align*}\left\| u_i^H(z)  - u_i^B(z) \right\|_{H_{x}^{s}L_z^2} \leq & \frac{C_1}{(N/2+1)^{q/2}} \\[4pt]
&  + \min\left\{ C_2 (\sqrt{N}+1)\, \e^{1-\alpha} \max \{\e^{\min(\ell, 1/2)}, V_{T}(\varepsilon)\}, \, C \sqrt{N}\delta \right\}. 
\end{align*}
If additionally $\delta=o(1/\sqrt{N})$, the above bound simplifies to 
\begin{align}
\begin{split}\label{Result2}\left\| u_i^H(z)  - u_i^B(z) \right\|_{H_{x}^{s}L_z^2} \leq & \frac{C_1}{(N/2+1)^{q/2}} \\[4pt]
&  + \min\left\{ C_2 (\sqrt{N}+1)\, \e^{1-\alpha} \max \{\e^{\min(\ell, 1/2)}, V_{T}(\varepsilon)\}, \, C_3 \right\}, 
\end{split}
\end{align}
with $C_3=o(1)$. 
\end{theorem}

\textbf{Proof. } First, we generalize the result  for the deterministic problem \cite[Theorem 8.1]{MB15} to the equation with uncertainty, and extend from the INS scaling to the acoustic scaling. If for all $z\in I_z$, the initial distribution $h_{\text{in}}\in \text{Null}(\mathcal L)$, $h_{\text{in}}\in H_{x,v}^{s} L_z^{\infty}$ and is sufficiently small, 
with some initial layer conditions, then for each $z$,  the function $h_{\e}$ converges strongly to 
\begin{equation}\label{h-eqn} h(t,x,v,z) = \left[\rho(t,x,z) + v\cdot u(t,x,z) + \frac{1}{2}(|v|^2 - d_v) \theta(t,x,z)\right] M(v) \end{equation}
in $L_{[0,T]}^2 H_x^s L_v^2$ as $\e\to 0$, where $\rho$, $u$, $\theta$ are the limit of the macroscopic quantities (density, mean velocity, temperature) obtained from $h_{\e}$.

If $h_{\text{in}}$ belongs to $H_x^{s+\ell}L_v^2 L_z^{\infty}$ for $\ell \in (0,1]$, one has the convergence estimate
\begin{equation}\label{Sup-h} \sup _{t \in[0, T]}\left\|h(t) -h_{\varepsilon}(t)\right\|_{H_{x}^{s} L_{v}^{2} L_z^{\infty}} \leqslant C\, \e^{1-\alpha} \max \left\{\e^{\min(\ell, 1/2)}, V_{T}(\varepsilon)\right\}, \end{equation}
where $T>0$ is the finite time. For each fixed $z$, the term $V_T(\e)$ is defined by
\begin{equation}\label{V_T}  V_T(\e)=\sup_{t\in [0, T]}||h(t) - h_{\e}(t)||_{L_x^{\infty}L_v^2} \to 0, \quad\text{ as   }\e\to 0. \end{equation}

The linearized collision operator $\mathcal L$ acts on $L_v^2$ with $\text{Ker}(\mathcal L) = \text{Span}\{\phi_1, \cdots, \phi_{d_v}\}$ where $(\phi_i)_{1\leq i \leq N}$ is an orthonormal family, $\pi_{\mathcal L}$ is the orthogonal projection on $\text{Ker}(\mathcal L)$ in $L_v^2$: 
\begin{equation}\label{Pi-L}\forall g \in L_{v}^{2}, \quad \pi_{\mathcal L}(g)=\sum_{i=1}^{N}\left(\int_{\mathbb{R}^{d_v}} g \phi_{i} d v\right) \phi_{i}, \end{equation}
where $\phi_{i}=P_{i}(v) e^{-|v|^{2}/4}$ and $P_i$ are polynomials with $P_0=1$, $P_1=v$, $P_2=\frac{1}{d_v}(|v|^2 - d_v)$. 
From \cite[equation (3.2)]{MB15}, by the Cauchy-Schwarz inequality, we know that 
$$ \forall s\in \mathbb N, \exists\, C>0, \forall h\in H_{x,v}^{s}, \qquad \left\|\pi_{\mathcal L}(h)\right\|_{H_{x, v}^{s}}^{2} \leqslant C\, \|h\|_{H_{x, v}^{s}}^{2}. $$
Take $\pi_{\mathcal L}$ on $h - h_{\e}$ in \eqref{Sup-h}, where $h$ is given by \eqref{h-eqn} which lies in Ker($\mathcal L$) and satisfies $\pi_{\mathcal L}(h)=h$, thus at fixed time $T>0$, 
$$ \left\| \pi_{\mathcal L}(h_{\e})(T) - \pi_{\mathcal L}(h)(T)\right\|_{H_{x}^{s} L_{v}^{2} L_z^{\infty}} \leq \left\|h(T) - h_{\varepsilon}(T)\right\|_{H_{x}^{s} L_{v}^{2} L_z^{\infty}} \leqslant C\, \e^{1-\alpha} \max \left\{\e^{\min(\ell, 1/2)}, V_{T}(\varepsilon)\right\}. $$
Here all the $C$ are generic constants. 
The macroscopic quantities are components of the summation in the formula \eqref{Pi-L}, therefore, at time $T>0$, 
$$ || u_i^H(T) - u_i^L(T) ||_{H_{x}^{s}L_z^{\infty}} \leq C\, \e^{1-\alpha} \max \left\{\e^{\min(\ell, 1/2)}, V_{T}(\varepsilon)\right\}. $$

Add up the projection error of the greedy algorithm, in particular \eqref{d-N}, we conclude the result in Theorem \ref{Thm2}. \qquad\qquad\qquad\qquad\qquad\qquad\qquad
\qquad\qquad\qquad\qquad $\square$

\begin{remark}
We summarized our error analysis results in \eqref{Result1} and \eqref{Result2}. Under the perturbative setting \eqref{PS}, assume 
{\bf Assumption 1} for the random collision kernel and initial data, in addition to initial condition \eqref{H-IC} satisfied, 
we conclude that our error estimate for all cases of $\e$ is an {\it uniform}-in-$\e$ result. 
\end{remark}

\section{Other choices of the low-fidelity models}
\label{sec:3}

As stated in Remark \ref{Rmk}, for the Boltzmann equation with random parameters, we can also choose some types of moment closure system, for example \cite{Ruo13}, as the low-fidelity model. In that case the error between the high-fidelity and bi-fidelity approximation is composed of the error in moment system solutions, compared with 
obtaining moments from the distribution function solved by the Boltzmann equation, and the projection error evoked by the greedy algorithm shown in \eqref{NW}. 
{\color{blue}This idea is currently carried out and numerically studied by the third author. }

{\color{blue}We discuss another possible way of choosing low-fidelity models in this section.}
One can let the low-fidelity model the same kinetic equation as the high-fidelity model, while solve it on a coarser temporal and physical mesh, 
in particular, using a larger $\Delta t$, $\Delta x$, $\Delta v$ than that for the high-fidelity model, under the constraint of the CFL condition. 
Given that it is too challenging to conduct error estimates of a discretized scheme for the full nonlinear Boltzmann equation, we will study the linear transport and linear Boltzmann equations to illustrate our idea. 

\subsection{The linear transport equation}

Consider the linear transport equation in one-dimensional slab geometry:
\begin{equation}
\label{LTE}
\left\{
\begin{array}{ll}
\displaystyle \varepsilon \partial_{t} f+v \partial_{x} f =\frac{\sigma}{\varepsilon} \mathcal{L} f, \qquad \sigma(x, z) \geq \sigma_{\min }>0, \\[4pt]
\displaystyle   \mathcal{L} f(t, x, v, z) =\frac{1}{2} \int_{-1}^{1} f\left(t, x, v^{\prime}, z\right) \mathrm{d} v^{\prime}-f(t, x, v, z), 
 \end{array}
 \right. 
 \end{equation}
where $f$ is the probability density distribution of particles at
position $x\in \Omega\subset\mathbb R$, time $t$, and with $v\in(-1,1)$ the cosine of the angle
between the particle velocity and its position variable, $\sigma(x,z)$ is the scattering coefficient depending on a random parameter 
$z\in I_z$, a $d$-dimensional random variable with probability distribution $\pi(z)$ known in priori. The parameter $\e$ is the Knudsen number defined as the ratio of the mean free path over a typical length scale such as the size of the spatial domain. 
The initial condition is uncertain $f(0,x,v,z)=f_0(x,v,z)$. 
Consider the periodic boundary condition $f(t,0,v,z) = f(t,1,v,z)$, for the sake of analysis. 
\\[2pt]

\noindent{\bf The diffusion limit}\, 
Denote $[\phi]=\frac{1}{2} \int_{-1}^{1} \phi(v) \mathrm{d} v$. 
Let the high-fidelity solution $$ u^H :=\rho=[f]. $$ For each fixed $z$, the classical diffusion limit theory of linear transport equation \cite{Bardos} applies, as $\e\to 0$, $\rho$ converges to the following random diffusion equation:
\begin{equation}\label{diffusion}\partial_{t} \rho=\partial_{x}\left(\kappa(z) \partial_{x} \rho\right), \end{equation}
where the diffusion coefficient $\kappa(z)=\frac{1}{3} \sigma(z)^{-1}$. 

Define the $\Gamma$-norm in space and velocity:
$$\|f(t, \cdot, \cdot, \cdot)\|_{\Gamma}^{2} := \int_{\Omega}\int_{-1}^1  \|f(t, x, v, \cdot)\|_{\pi}^{2}\,\mathrm{d} v \mathrm{d} x, \quad t \geq 0. $$
In \cite{MM08}, Lemou and Mieussens have proposed an AP scheme for the linear kinetic equations in the diffusion limit, 
based on the micro-macro decomposition of the distribution function into its microscopic and macroscopic components. 

By the micro-macro decomposition method, one lets
$$ f = \rho + \e g, $$
where $\rho = [f]$ and $g$ is such that $[g]=0$, then the micro-macro model reads
\begin{equation}
\label{mm}
\left\{
\begin{array}{ll}
\displaystyle \partial_t \rho + \partial_x [vg] = - \sigma \rho, \\[4pt]
\displaystyle \partial_t g + \frac{1}{\e}(I - [\cdot])(v\partial_x g) = - \frac{\sigma}{\e^2}Lg - \frac{1}{\e^2}v \partial_x \rho, 
\end{array}
\right. 
\end{equation}
with initial data $\rho_{in} = [f_{in}]$ and $\e g_{in} = f_{in} - \rho_{in}$. 

\cite[Section 4.2]{JGL10} studied the analysis of the convergence error in the above micro-macro decomposition framework. 
The study of uncertainty was done in \cite{Jin-Liu-Ma}. 
Define the convergence errors for fixed $z$ in $\rho$ and $g$, and fixed $v$ in $g$, 
$$\tilde\rho_i^n = \rho(t_n, x_i, z) - \rho_i^n(z), \qquad
\tilde g_{i+\frac{1}{2}}^n(v,z) = g(t_n, x_{i+\frac{1}{2}},v,z) - g_{i+\frac{1}{2}}^n(v,z). $$
We will omit the $v$ dependence in $g$. $\rho(t_n, x_i)$, $g(t_n, x_{i+\frac{1}{2}})$ are the true solutions to the micro-macro system \eqref{mm} 
at $(t_n, x_i)$ and $(t_n, x_{i+\frac{1}{2}})$ respectively, $\rho_i^n$, $g_{i+\frac{1}{2}}^n$ are the corresponding numerical solutions. 

Define for grid functions $\mu = (\mu_i)_{i\in\mathbb Z}$ and 
velocity dependent grid function $\phi(v)=(\phi_{i+\frac{1}{2}}(v))_{i\in\mathbb Z}$, 
\begin{equation} || \mu ||^2 = \sum_{i\in \mathbb Z}\mu_i^2\,\Delta x, \qquad
 \left| || \phi || \right| = \sum_{i \in\mathbb Z} \left[ \phi_{i+\frac{1}{2}}^2\right] \Delta x. \end{equation}
The authors in \cite{JGL10} showed that
\begin{equation}\label{Energy} || \tilde \rho^n || + \e \left| || \tilde g^n || \right| 
\leq C \left( (1+\e^2)\Delta t + \Delta x^2 + \e \Delta x\right), \end{equation}
where $C$ is independent of $\Delta t$, $ \Delta x$ and $\e$. 

In our bi-fidelity method, the numerical solutions $\rho^H$, $g^H$ to the discretized system \eqref{mm} are considered as the high-fidelity solutions, with the 
discretization parameters $\Delta x_1$, $\Delta t_1$. 
The numerical solution $\rho^L$, $g^L$ to \eqref{mm} are the low-fidelity solutions, solved by using coarser mesh in the physical space $\Delta x_2$, $\Delta t_2$, 
with $\Delta x_2 > \Delta x_1$, $\Delta t_2 > \Delta t_1$. 
Adapt \eqref{Energy} to the low- and high-fidelity solutions at time $t_n$, 
\begin{align*}
& || \tilde \rho^{H,n} || + \e \left| || \tilde g^{H,n} || \right| 
\leq C_1 \left( (1+\e^2)\Delta t_1 + \Delta x_1^2 + \e \Delta x_1\right), \\[4pt]
& || \tilde \rho^{L,n} || + \e \left| || \tilde g^{L,n} || \right| 
\leq C_2 \left( (1+\e^2)\Delta t_2 + \Delta x_2^2 + \e \Delta x_2\right). 
\end{align*}

Since at each $(t_n, x_i)$, $$ \rho_i^{H,n}(z) - \rho_i^{L,n}(z) = \rho_i^{H,n}(z) - \rho(t_n, x_i, z) - \left( \rho_i^{L,n}(z) - \rho(t_n, x_i, z)\right), $$
and similarly for the difference between $g_{i+\frac{1}{2}}^{H,n}$ and $g_{i+\frac{1}{2}}^{L,n}$, then by the Cauchy-Schwarz inequality 
for each fixed $z$, we get 
\begin{align*}
&\quad \left\| \rho^{H,n}(z) - \rho^{L,n}(z) \right\| + \e\, \left| \left\| g^{H,n}(z) - g^{L,n}(z)\right\| \right|  \\[4pt]
& \leq C_0 \left( (1+\e^2)\Delta t_2 + \Delta x_2^2 + \e \Delta x_2\right) + ``\text{interpolation errors}" \\[4pt]
& \leq C\left( (1+\e^2)\Delta t_2 + \Delta x_2^2 + \e \Delta x_2\right), 
\end{align*}
where $C_0 = 2 \max\{ C_1, C_2\}$, and the small interpolation errors between different meshes is due to that the low- and high-fidelity models may not be solved simultaneously at the same time $t_n$ and spatial point $x_i$, however this error can be bounded by the larger discretization parameters. 

The analyticity of the solution to \eqref{LTE} with respect to the uncertainty is proved in 
\cite{Qin-Li}, under suitable assumptions for the random scattering coefficient and initial data. We paraphrase it here: 
\begin{proposition}\cite[Theorem 8]{Qin-Li}
If we assume that 
$$ |\partial_z \sigma(x,z)| \leq C_1, \qquad \partial_z^l \sigma(x,z)=0 \, \text{  for  } l\geq 2. $$
Denote $g_l = \partial_z^l f$, if the initial data satisfies \begin{equation}\label{g-IC} || g_l(t=0, z)|| \leq H^l, \qquad \text{for all  } l \geq 0, \end{equation}
then $g_l$ at time $t$ has the following estimate for all $l\geq 0$, 
$$ || g_l || \leq C e^{-\lambda t} \left( H + t \frac{\tilde C_1}{\e^2}\right)^l. $$
\end{proposition}
The above result indicates that the Taylor series $$f(z)=\sum_{l=0}^{\infty} \frac{g_{l}}{l !}\left(z-z_{0}\right)^{l}$$
converges for any $\e$, with the convergence radius $$ r\left(z_{0}\right)=\frac{1}{\limsup _{l \rightarrow \infty}\left(g_{l}\left(z_{0}\right) / l !\right)^{1 / l}} = \infty, $$ 
thus the solution $f$ is analytic (or holomorphic) with respect to the randomness. This analysis is performed in an one-dimensional fashion for each direction $z_n$, in general, we assume that $\sigma(x,z)$ has an {\it affine dependence} on $z$: 
\begin{equation}\label{Sigma} \sigma(x,z) = \sigma_0(x) + \sum_{j \geq 1} z_j \psi_j(x),  \qquad z:=(z_j)_{j\geq 1}, \end{equation}
where the sequence $\left(\left\|\psi_{j}\right\|_{L^{\infty}(x)}\right)_{j \geq 1} \in \ell^{p}$ for $0<p<1$. 
With this holomorphicity, one can adopt the theory in \cite{Cohen15}, and in particular 
\eqref{NW} and \eqref{d-N} as in Section \ref{sec:2}. Therefore, the overall error of the bi-fidelity method with choosing coarse/fine meshes in the low- and high-fidelity models and using the micro-macro decomposition numerical scheme is given on a discretized level (in $t$ and $x$) by
\begin{align*}
&\quad \left\| \rho^{H,n}(z) - \rho^{B,n}(z) \right\|_{L^2_z} + \e\, \left| \left\| g^{H,n}(z) - g^{B,n}(z)\right\| \right|_{L^2_z}  \\[4pt]
& \leq  \frac{C}{(N/2+1)^{q/2}} + C^{\prime}\, (\sqrt{N}+1) \left( (1+\e^2)\Delta t_2 + \Delta x_2^2 + \e \Delta x_2\right), 
\end{align*}
where $C$, $C^{\prime}$ are independent of $\e$, and $q$ is associated to the affine representer $\psi_j$ in \eqref{Sigma}, with the relation $q=\frac{1}{p}-1$.  
This error estimate result is uniform in $\e$. 

\begin{theorem}
If the initial data $\partial_z^l f$ satisfies \eqref{g-IC}, the random scattering coefficient is affine in $z$ as shown in \eqref{Sigma}, 
then the error between the high-fidelity and bi-fidelity solutions solved by the micro-macro system \eqref{mm} has the following estimate: 
\begin{align}
\begin{split}
&\quad \left\| \rho^{H}(z) - \rho^{B}(z) \right\|_{L^2_z} + \e \left| \left\| g^{H}(z) - g^{B}(z)\right\| \right|_{L^2_z}  \\[4pt]
& \leq  \frac{C}{(N/2+1)^{q/2}} + C^{\prime}\, (\sqrt{N}+1) \left( (1+\e^2)\Delta t_2 + \Delta x_2^2 + \e \Delta x_2\right), 
\end{split}
\label{Error-LTE}
\end{align}
where $\Delta x_2$, $\Delta t_2$ are the coarser numerical meshes in obtaining the low-fidelity solutions, 
while $\Delta x_1$, $\Delta t_1$ are used in the high-fidelity solvers. 
\end{theorem}

\subsection{The linear Boltzmann equation with anisotropic scattering kernel}

Consider the linear Boltzmann equation with uncertain anisotropic scattering kernel and the external potential $V(x)$ that is independent of $z$, 
\begin{equation}\label{LBE}\frac{\partial f}{\partial t} + \alpha \cdot \nabla f = Q(f), \end{equation}
where
$$Q_{\sigma}(f)(t,x,v,z) =\int_{v^{\prime} \in \mathcal{D}_{v}^{d}}\left(\sigma\left(x, v, v^{\prime}, z\right) f^{\prime}-\sigma\left(x, v^{\prime}, v, z\right) f\right) d v^{\prime}. $$
where $f^{\prime}=f(x,v^{\prime},z)$, and \(\begin{array}{lll} {\alpha(x, v)=\left(\begin{array}{c}{v} \\ {-E(t, x)}\end{array}\right),} & {E=-\nabla_{x} V, \quad \nabla=\left(\begin{array}{c}{\nabla_{x}} \\ {\nabla_{v}} \end{array}\right)}\end{array}. \)
It has been shown in \cite{Qin-Li} that the solution to \eqref{LBE} is analytic in the random space, if $\sigma$ has either an affine dependence on $z$ or arbitrary dependence on $z$ with $\left|\frac{\partial_z^n \sigma}{n!}\right|$ uniformly bounded. 

We recall the $L^2$ error estimate of the semi-discrete discontinuous Galerkin (DG) discretization of \eqref{LBE} shown in \cite[Theorem 9]{Cheng-Gamba}, in a special case when $M(v)=\text{constant}$: 
$$ \| f_{h}(t, \cdot) -f(t, \cdot, \cdot) \|_{L^2\left(\Omega_{\mathcal{D}}\right)} \leq
 C \sqrt{t}\, e^{C h t}\, h^{k+\frac{1}{2}} |f|_{L^{\infty}\left([0, t], H^{k+1}\left(\Omega_{\mathcal{D}}\right)\right)}, 
$$
where $C=C\left(\operatorname{diam}\left(\Omega_{\mathcal{D}}\right),\|\alpha\|_{\mathbf{W}^{\{1, \infty\}}\left(\Omega_{\mathcal{D}}\right)}\right)$ with $\Omega_{\mathcal{D}}$ the computational domain, and $C$ is independent of $h$ or $t$. 

Similar to \eqref{Error-LTE}, we obtain the error estimate of the bi-fidelity method for the linear Boltzmann equation \eqref{LBE} when choosing the low-fidelity model as the same equation computed on coarser mesh $\Delta x_2$, 
\begin{align*}
\begin{split}
 \| f^H - f^B \|_{L^2\left(\Omega_{\mathcal{D}}\right) L^2_z} & \leq \frac{C}{(N/2+1)^{q/2}}   \\[4pt]
&\quad + C^{\prime}(\sqrt{N}+1) \sqrt{t}\, e^{c \Delta x_2 t}\, ({\Delta x_2})^{k+\frac{1}{2}} |f|_{L^{\infty}\left([0, t], H^{k+1}\left(\Omega_{\mathcal{D}}\right)\right)}. 
\end{split}
\end{align*}

\section{Conclusion}
\label{sec:4}
In this paper, we provide a general framework to 
study error estimates of the bi-fidelity method to solve multi-scale kinetic equations with random inputs. Our result is {\it uniform} in the Knudsen number.  
We take the Boltzmann and the linear transport equations as 
two important classes of examples to conduct the analysis. Our idea and analysis can be applied to general multiscale kinetic equations with uncertainty since the 
regularity analysis in the previous work \cite{LJ-UQ} works for general linear or non-linear kinetic equations. 
Our work sheds some lights on choosing low-fidelity model for the bi-fidelity method to solve kinetic equations with random parameters. 

\begin{appendices}

\section{Proof of \eqref{h-eps}}
We first discuss the incompressible Navier-Stokes (INS) scaling. The kinetic equation reads
$$\partial_t f + \frac{1}{\varepsilon} v\cdot\nabla_x f = \frac{1}{\varepsilon^2} \mathcal Q(f, f), $$
and the perturbative setting
\begin{equation}\label{pert} f = M + \delta \sqrt{M} h, \end{equation}
then $h$ satisfies
$$ \partial_t h + \frac{1}{\varepsilon} v\cdot\nabla_x h = \frac{1}{\varepsilon^2}\mathcal L(h) + \frac{\delta}{\varepsilon^2} \mathcal F(h, h). $$
Following a similar proof as \cite{MB15} for the deterministic problem and refer to \cite{LJ-UQ} with the randomness, while consider the setting \eqref{pert} with the small parameter $\delta$ 
(in \cite{LJ-UQ} $\delta = \varepsilon$), one finally gets
$$\frac{d}{d t}\|h\|_{\mathcal{H}_{\varepsilon \perp}^{s,r}}^{2} \leqslant\left( \frac{\delta^2}{\varepsilon^2}\, K_1 C_{\Gamma}^{2} C\|h\|_{\mathcal{H}_{\varepsilon \perp}^{s,r}}^{2}-K_0\right)\|h\|_{H_{\Lambda}^{s,r}}^{2}. $$
If 
$$ || h_{\text{in}} ||_{\mathcal{H}_{\varepsilon \perp}^{s,r}}^2 \leq \frac{\varepsilon^2}{\delta^2}\frac{K_0}{2 K_1 C_{\Gamma}^{2} C }, $$
then $ \|h\|_{\mathcal{H}_{\varepsilon \perp}^{s,r}}$ always decreases, and for all $t>0$, 
$$\frac{d}{d t}\|h\|_{\mathcal{H}_{\varepsilon \perp}^{s,r}}^{2} \leqslant-\frac{K_0}{2}\|h\|_{H_{\Lambda}^{s,r}}^{2}, $$
since the $H_{\Lambda}^{s,r}$-norm controls the $H_{x,v}^{s,r}$-norm which is equivalent to the 
$\mathcal{H}_{\varepsilon \perp}^{s,r}$-norm, applying Gronwall's lemma gives us the exponential decay in the $\mathcal{H}_{\varepsilon \perp}^{s,r}$-norm. 

To conclude, if \begin{equation}\label{h-in} || h_{\text{in}} ||_{\mathcal{H}_{\varepsilon \perp}^{s,r}} \leq \frac{\varepsilon}{\delta}\,\eta, \end{equation}
then for all $t>0$, \begin{equation}\label{h-T}\|h\|_{\mathcal{H}_{\varepsilon \perp}^{s,r}} \leq \eta^{\prime}\, e^{-\tau t}, \end{equation}
where $\eta$ is sufficiently small, and the constants $\eta$, $\eta^{\ast}$, $\tau$ are all independent of $\varepsilon$ and $\delta$. 
We don't specify the $r$-dependence of $\eta$ and $\eta^{\prime}$ though. 
Since \(\|\cdot\|_{\mathcal{H}_{\varepsilon \perp}^{s,r}} \sim\|\cdot\|_{H_{x, v}^{s,r}}\), then \eqref{h-in}--\eqref{h-T} also hold true for $||h||_{H_{x,v}^{s,r}}$. 

To conclude, consider the general kinetic problem with uncertainty and scalings, 
\begin{equation*}
\left\{
\begin{array}{ll}
\displaystyle\partial_{t} f+\frac{1}{\varepsilon^{\alpha}} v \cdot \nabla_{x} f=\frac{1}{\varepsilon^{1+\alpha}} \mathcal{Q}(f, f), \\[4pt]
\displaystyle f(0, x, v, z)=f_{\text{in}}(x, v, z), \quad x \in \Omega \subset \mathbb{T}^{d}, v \in \mathbb{R}^{d}, z \in I_{z}  \subset \mathbb R, 
 \end{array}
 \right.
 \end{equation*}
where $\alpha=0$ and $\alpha=1$ correspond to the acoustic and INS scalings, respectively. 

If we assume the initial data 
$$|| h_{\text{in}} ||_{H_{x,v}^s L_z^{\infty}} \leq \frac{\varepsilon}{\delta}\,\eta,$$
with sufficiently small $\eta_s$, {\color{red}and the random collision kernel satisfy \eqref{Coll},} then at all time $t>0$, 
\begin{itemize}
 \item (i) For $\alpha=1$ (INS scaling), 
$$\left\|h_{\varepsilon}\right\|_{H_{x, v}^{s, r} L_{z}^{\infty}} \leq \eta^{\prime}\, e^{-\tau t}, \quad\left\|h_{\varepsilon}\right\|_{H_{x, v}^{s} H_{z}^{r}} \leq \eta^{\prime}\, e^{-\tau t}, $$
\item (ii) For $\alpha=0$ (acoustic scaling), 
$$\left\|h_{\varepsilon}\right\|_{H_{x, v}^{s, r} L_{z}^{\infty}} \leq \eta^{\prime}\, e^{-\varepsilon \tau t}, \quad\left\|h_{\varepsilon}\right\|_{H_{x, v}^{s} H_{z}^{r}} \leq \eta^{\prime}\, e^{-\varepsilon \tau t}, $$
\end{itemize}
 where all the constants $\eta$, $\eta^{\prime}$, $\tau$ are independent of $\varepsilon$ and $\delta$.

\section{Proof of analyticity for the linearized Boltzmann equation}

We prove the analyticity in the random space for the linearized Boltzmann equation with multiple scales. 
The perturbative solution satisfies 
\begin{equation}\label{h-App} \partial_t h + \frac{1}{\e^{\alpha}} v \cdot \nabla_x h = \frac{1}{\e^{1+\alpha}} \mathcal L(h), \end{equation}
where $\alpha=0$ corresponds to the acoustic scaling, and $\alpha=1$ corresponds to the INS scaling. 

For simplicity, we provide a one-dimensional analysis for each direction of $z_n$. 
Assume the collision cross-section $b$ has an affine dependence on $z$, 
$$ b(\cos\theta) = b_0(\cos\theta) + b_1(\cos\theta)\, z, \qquad |b_1|\leq C_b. $$ 
Define the operator $\mathcal T_{\e}:= \frac{1}{\e^{1+\alpha}}\mathcal L - \frac{1}{\e^{\alpha}} v\cdot\nabla_x$, and denote
$h_l := \partial_z^l h$ for all $l\geq 0$. 
Take $\partial_z^l$ of the equation \eqref{h-App}, then 
$$ \partial_t h_l = \mathcal T_{\e}(h_l) + \frac{1}{\e^{1+\alpha}} \mathcal L_{b_1}(h_{l-1}), $$
where $\mathcal L_{b_1}$ is defined by substituting $b(\cos\theta)$ by $b_1(\cos\theta)$ in the linearized collision operator $\mathcal L$. 
Thus
$$ \frac{d}{dt} ||h_l||_{H^k}^2 = 2 \langle \mathcal T_{\e}(h_l), h_l \rangle_{H^k} + \frac{2}{\e^{1+\alpha}} l \, \langle \mathcal L_{b_1}(h_{l-1}), h_l \rangle_{H^k}. $$
By the hypocoercivity theory and $\mathcal L$ being a bounded operator, 
$$  \frac{d}{dt} ||h_l||_{H^k}^2  \leq - \e C_0\, ||h_l||_{H_{\Lambda}^k}^2 + \frac{C_1}{\e^{1+\alpha}} l \, ||h_{l-1}||_{H^k} ||h_l||_{H^k}, $$
where 
$$\|h\|_{H_{\Lambda}^k}^2=\sum_{|j|+|l|=k}\left\|\frac{\partial^2}{\partial v_j \partial x_l}h\right\|_{\Lambda}^2, $$
and $$\|h\|_{\Lambda}=\left\|h(1+|v|)^{\gamma / 2}\right\|_{L^{2}_{x,v}}, $$ with $\gamma$ shown in \eqref{Phi} on the 
potential part of the collision kernel. 
Set $g_l = ||h_l||_{H^k}$, since $ ||h_l||_{H_{\Lambda}^k}$ controls $||h_l||_{H^k}$, then 
$$ \frac{d}{dt}g_l \leq -\e C_0\, g_l + \frac{C_1}{\e^{1+\alpha}} l\, g_{l-1}. $$
Using \cite[Lemma 6]{Qin-Li} in a similar way, one gets
$$ g_l(t) \leq e^{-\e C_0 t} \sum_{m=0}^l \frac{l!}{(l-m)!\, m!}\left(\frac{C_1}{\e^{1+\alpha}}t\right)^m\, h_{l-m}(0). $$
 Following \cite[Theorem 5]{Qin-Li}, if the initial data satisfies 
 $ ||\partial_z^l h_0(z) ||_{H^k} \leq R^l$, $\forall l \geq 0$, 
 then at time $t$, 
 $$ ||\partial_z^l h(z)||_{H^k} \leq C e^{-\e C_0 t}\left( R + \frac{C_1}{\e^{1+\alpha}} t \right)^l, $$
 where $C$, $C_0$, $C_1$ are independent of $l$ and $\e$. 
  The convergence radius for $h$ at any point $z_0$ is defined by
 $$r\left(z_{0}\right)=\frac{1}{\limsup _{l \rightarrow \infty}\left(g_{l}(z_{0})/\, l !\right)^{1 / l}}=\infty, $$
 which is independent of $z_0$, thus $h$ (or $f$) is {\it analytic} in $z$. 
 
\end{appendices}

\bibliographystyle{siam}
\bibliography{MultiFid2.bib}

\end{document}